\documentclass[runningheads]{llncs}

\usepackage{graphicx,color}
\graphicspath{ {./figure/} }

\usepackage{amsmath,amssymb,latexsym}
\usepackage{mathtools}
\usepackage{xspace}
\usepackage{tikz}
\usepackage{subcaption}
\usepackage[utf8]{inputenc}
\usepackage{hyperref}

\newcommand{\BF}{\mathit{BF}}
\newcommand{\CCC}{\mathit{CCC}}

\newcommand{\cp}{\,\square\,}

\newcommand{\bp}{\mathrm{bp}}
\newcommand{\Prop}{\mathcal{P}}

\newcommand{\comment}[2]{ { \noindent \color{red}{\small [$\bullet$ \textsc{#1}: \textsl{#2}]}}}

\begin{document}
\title{Mutual visibility in hypercube-like graphs\thanks{The work has been supported in part by the Italian National Group for Scientific Computation (GNCS-INdAM).}
}
%
%
\author{%
Serafino Cicerone\inst{1}
\and
Alessia {Di Fonso}\inst{1}
\and
Gabriele {Di Stefano}\inst{1}
\and
Alfredo Navarra\inst{2}
\and
Francesco Piselli\inst{2}
}
\authorrunning{S. Cicerone et al.}
%

\institute{Dipartimento di Ingegneria e Scienze dell'Informazione e Matematica,\\
        Università degli Studi dell'Aquila, I-67100 
        L'Aquila, Italy. \\
\email{serafino.cicerone@univaq.it}, \email{alessia.difonso@univaq.it}, \email{gabriele.distefano@univaq.it}
\and 
Dipartimento di Matematica e Informatica,
        Università degli Studi di Perugia I-06123 
        Perugia, Italy.
\email{alfredo.navarra@unipg.it},
\email{francesco.piselli@unifi.it}
}

\maketitle              

\sloppy 

\begin{abstract}
Let $G$ be a graph and $X\subseteq V(G)$. Then, vertices $x$ and $y$ of $G$ are $X$-visible if there exists a shortest $u,v$-path where no internal vertices belong to $X$. The set $X$ is a mutual-visibility set of $G$ if every two vertices of $X$ are $X$-visible, while $X$ is a total mutual-visibility set if any two vertices from $V(G)$ are $X$-visible. The cardinality of a largest mutual-visibility set (resp. total mutual-visibility set) is the mutual-visibility number (resp. total mutual-visibility number) $\mu(G)$ (resp. $\mu_t(G)$) of $G$. It is known that computing $\mu(G)$  is an NP-complete problem, as well as $\mu_t(G)$. In this paper, we study the (total) mutual-visibility in hypercube-like networks (namely, hypercubes, cube-connected cycles, and butterflies). Concerning computing $\mu(G)$,
we provide approximation algorithms for both hypercubes and cube-connected cycles, while we give an exact formula for butterflies. Concerning computing $\mu_t(G)$ (in the literature, already studied in hypercubes), we provide exact formulae for both cube-connected cycles and butterflies.   

\end{abstract}
%
%
\section{Introduction}\label{sec:introduction}
%
Problems about sets of points in the Euclidean plane and their mutual visibility have been investigated for a long time. For example, in~\cite{Dudeney17} Dudeney posed the famous \emph{no-three-in-line} problem: finding the maximum number of points that can be placed in an $n \times n$ grid such that there are no three points on a line. Beyond the theoretical interest, solutions to these types of geometric/combinatorial problems have proved useful in recent times in the context of swarm robotics. The requirement is to define algorithms that allow autonomous mobile robots to change in a finite time their configuration in the plan so they can see each other (see, e.g.,~\cite{LunaFCPSV17,PoudelAS21,SharmaVT21}).

Mutual visibility in graphs with respect to a set of vertices has been recently introduced and studied in~\cite{DiStefano22} in the sense of the existence of a shortest path between two vertices not containing a third vertex from such a set. The visibility property is then understood as a kind of non-existence of ``obstacles'' between the two vertices in the mentioned shortest path, which makes them ``visible'' to each other. 
Formally, let $G$ be a connected and undirected graph, and $X\subseteq V(G)$ a subset of the vertices of $G$. Two vertices $x, y \in V(G)$ are \emph{$X$-visible} if there exists a shortest $x,y$-path where no internal vertex belongs to $X$. $X$ is a \emph{mutual-visibility set} if its vertices are pairwise $X$-visible. The cardinality of a largest mutual-visibility set is the \emph{mutual-visibility number} of $G$, and it is denoted by $\mu(G)$. Computing one of such largest sets, referred to as a \emph{$\mu$-set} of $G$, solves the so-called \textsc{Mutual-Visibility} problem.  In~\cite{DiStefano22}, it is also shown that computing $\mu(G)$ is an NP-complete problem, but there exist exact formulae for the mutual-visibility number of special graph classes like paths, cycles, blocks, cographs, and grids. 
In the subsequent works~\cite{CiceroneDK23} and~\cite{CiceroneDKY23}, exact formulae have been derived also for both the Cartesian and the Strong product of graphs. The contributions presented in those initial works showed several interesting connections with other mathematical contexts like the relationship existing between such a visibility problem and an instance of the very well-known Zarankiewicz problem (see \cite[Corollary 3.7]{CiceroneDK23}). The mutual-visibility is also related to the general position problem in graphs (see, e.g.,~\cite{ManuelK18,KlavzarNC22,Prabha23} and references therein). It is worth noting that this new definition of mutual visibility on graphs has already been applied in the context of swarm robotics where robots operate on environments modelled by graphs and where visibility is verified along shortest paths (cf.~\cite{CiceroneFSN23,CDDN2023-gmv-grids}). 

In~\cite{CiceroneDKY23}, a natural extension of the mutual-visibility has been proposed. Formally,  $X\subseteq V(G)$ is a \emph{total mutual-visibility set} of $G$ if every two  vertices $x$ and $y$ of $G$ are $X$-visible. A largest total mutual-visibility set of $G$ is a $\mu_t$-set of $G$, and its cardinality is the \emph{total mutual-visibility number} of $G$ denoted as $\mu_t(G)$. Of course, $\mu_t(G) \le \mu(G)$. In~\cite{variety23} it is shown that also computing $\mu_t(G)$ is an NP-complete problem. 

The most recent works in the literature on mutual-visibility concern total mutual-visibility. This setting is clearly more restrictive, but it surprisingly turns out to become very useful when considering networks having some Cartesian properties in the vertex set, namely, those ones of product-like structures (cf.~\cite{CiceroneDKY23}). 
In~\cite{totalmutualzero}, in fact, the total mutual-visibility number of Cartesian products is bounded and several exact results proved. Furthermore, a sufficient and necessary condition is provided for asserting when $\mu_t(G) = 0$.
In~\cite{KuziakR23}, the authors give several bounds for $\mu_t(G)$ in terms of the diameter, order and/or connected domination number of $G$. They also determine the exact value of the total mutual-visibility number of lexicographic products. Finally, the total mutual-visibility of Hamming graphs (and, as a byproduct, for hypercubes) is studied in~\cite{BujtaKT23}.

\smallskip\noindent
\textbf{Results.}
In this work, we study the mutual-visibility and the total mutual-visibility of hypercube-like graphs. These graphs usually model interconnection networks which are widely used in parallel and distributed computing systems. They offer efficient communication and connectivity patterns that make them suitable for various applications, including parallel processing, supercomputing, and designing efficient interconnection topologies for computer systems. 
In these systems, two nodes that are $X$-visible can communicate in an efficient way, that is through shortest paths, and their messages can be maintained confidential: the exchanged messages do not pass through nodes in $X$.

A \emph{hypercube} of dimension $d$ is denoted as $Q_d$ and consists of $2^d$ vertices, each labelled by a binary string with $d$ bits. Two vertices are adjacent if and only if the two binary strings differ by exactly one bit.
A \emph{cube-connected cycle} of dimension $d$ is denoted as $\CCC_d$ and consists of $d\cdot 2^d$ vertices. It models networks that combine the properties of both hypercube graphs and cycle graphs. Informally, $\CCC_d$ can be obtained by $Q_d$ by replacing each vertex by a cycle with $d$ vertices. 
A \emph{butterfly} of dimension $d$ is denoted as $\BF(d)$ and consists of $(d+1)\cdot 2^d$ vertices. It is composed of interconnected stages of switches that resemble the wings of a butterfly. Each stage connects nodes in a structured way, allowing efficient communication between nodes at different levels. Butterfly networks are known for their logarithmic diameter and suitability for sorting and data exchange operations.
For these classes of graphs, we provide the following results:
\begin{itemize}
\item 
For hypercubes, we first prove the bounds $\mu(Q_d)\ge 
\binom{d}{\lfloor\frac d 2 \rfloor} + 
\binom{d}{\lfloor\frac d 2 \rfloor+3}$ and $\mu(Q_d)\le 2^{d-1}$, and then we exploit such results to get an $O(\sqrt{d})$-approximation algorithm for computing $\mu(Q_d)$. Alternatively, for an $n$-vertex hypercube, the approximation ratio can be expressed as $O( \sqrt{\log n} )$.
\item 
For cube-connected cycles, the provided bounds are $\mu(\CCC_d) \le 3\cdot 2^{d-2}$ and $\mu(CCC_d)\ge 
2^{\lceil \frac d 2 \rceil-1}$. They give rise to a $3 \cdot 2^{\lfloor \frac d 2 \rfloor-1}$-approximation algorithm (whose ratio can be expressed as $O(\sqrt{n})$ for an $n$-vertex cube-connected cycle). Concerning the total mutual-visibility, we prove that $\mu_t(\CCC_d)=0$.
\item 
For butterflies, we are able to provide exact formulae: $\mu( \BF(d) ) = 2^{d+1}-2$ and $\mu_t( \BF(d) ) = 2^d$.
\end{itemize}

\noindent \textbf{Outline.} 
Concerning the organization of the paper, in the next section we provide all the necessary notation and preliminary concepts. The subsequent Sections~\ref{sec:hypercubes},~\ref{sec:ccc}, and~\ref{sec:butterflies} are specialized for presenting - in order - results about (total) mutual visibility for hypercubes, cube-connected cycles, and butterflies. Concluding remarks are provided in Section~\ref{sec:conclusions}.

\section{Notation and preliminaries}\label{sec:notation}
%
In this work, we consider undirected and connected graphs. We use standard terminologies from~\cite{BrandstadtLS99}, some of which are briefly reviewed here.

Given a graph $G$, $V(G)$ and $E(G)$ are used to denote its vertex set and its edge set, respectively, and $n(G)$ is used to represent the size of $V(G)$. If $u,v\in V(G)$ are adjacent, $(u,v)\in E(G)$ represents the corresponding edge. 
If $X\subseteq V(G)$, then $G[X]$ denotes the subgraph of $G$ \emph{induced} by $X$, that is the maximal subgraph of $G$ with vertex set $X$. 

The usual notation for representing special graphs is adopted. $K_n$ is the \emph{complete graph} with $n$ vertices, $n\ge 1$, that is the graph where each pair of distinct vertices are adjacent. $P_n$ represents any \emph{path} $(v_1,v_2,\ldots,v_n)$ with $n\ge 2$ distinct vertices where $v_i$ is adjacent to $v_j$ if $|i-j|=1$. Vertices $v_2,\ldots,v_{n-1}$ are all \emph{internal} of $P_n$. A {\em cycle} $C_n$, $n\geq 3$, in $G$ is a path $(v_0,v_1,\ldots, v_{n-1})$ where also $(v_{0},v_{n-1})\in E(G)$. Two vertices $v_i$ and $v_j$ are {\em consecutive} in $C_n$ if  $j=(i+1)\bmod n$ or $i=(j+1)\bmod n$. The \emph{distance} function on a graph $G$ is the usual shortest path distance.


Two graphs $G$ and $H$ are \emph{isomorphic} if there exists a bijection $\varphi :V(G) \rightarrow V(H)$ such that $(u,v) \in E(G) \Leftrightarrow (\varphi(u),\varphi(v)) \in E(H)$ for all $u,v \in V(G)$. Such a bijection $\varphi$ is called \emph{isomorphism}.
Given $G$ and $H$, we consider also the following graph operation: the {\em Cartesian product} $G\cp H$ has the vertex set $V(G)\times V(H)$ and the edge set $E(G\cp H)  = \{((g,h),(g',h')):\ (g,g')\in E(G)\mbox{ and } h=h', \mbox{ or, } g=g' \mbox{ and }  (h,h')\in E(H)\}$.


%

If $H$ is a sugbraph of $G$, $H$ is said to be \emph{convex} if all shortest paths in $G$ between vertices of $H$ actually belong to $H$. Concerning convex subgraphs, we recall from~\cite{DiStefano22} the following useful statement.

\begin{lemma}{\rm \cite[Lemma~2.1]{DiStefano22}}
\label{lem:convex-subgraph}
Let $H$ be a convex subgraph of any graph $G$, and let $X$ be a mutual-visibility set of $G$. Then, $X\cap V(H)$ is a mutual-visibility set of $H$.
\end{lemma}

Finally, note that binary strings are used as vertex labels or components of vertex labels for the classes of graphs studied in this paper (hypercubes, cube-connected cycles, and butterflies). For a binary string $x = x_0x_1 \cdots x_i \cdots x_{d-1}$ with $d$ bits,  position 0 corresponds to the leftmost bit, and position $n-1$ to the rightmost bit. Sometimes, we interpret these strings as (binary) numbers. We use $x(i)$ to denote the binary string obtained from $x$ 
by complementing the bit in position $i$. We will abbreviate ``the vertex with label $x$'' to ``vertex $x$''. 

\section{Hypercube}\label{sec:hypercubes}
\begin{figure}[t]
\centering
\includegraphics[scale=0.60]{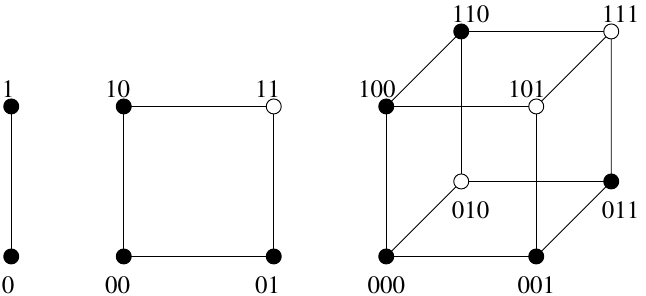}~~  
\caption{A representation of the hypercube $Q_d$ for $d = 1$, $2$, $3$. In each graph, the black vertices form a $\mu$-set.}
\label{fig:hypercube}
\end{figure}

%
 The $d$-dimensional hypercube $Q_d$ is an perperundirected graph with vertex set $V(Q_d) = \{0, 1\}^d$, and two vertices are adjacent if and only if the two binary strings differ by exactly one bit, that is, the Hamming distance%
 \footnote{The Hamming distance between two strings of equal length is the number of positions at which the corresponding symbols are different. In other words, it measures the minimum number of substitutions required to change one string into the other.}
 of the two binary strings is 1.
 It is worth noting that $Q_d$ can also be recursively defined in terms of the Cartesian product of two graphs as follows:
\begin{itemize}
\item
$Q_1 = K_1$,
\item
$Q_d= Q_{d-1} \cp K_2$, for $d\ge 2$.
\end{itemize}
This makes clear that $Q_d$ can be also seen as formed by two subgraphs both isomorphic to $Q_{d-1}$ (cf. Fig.~\ref{fig:hypercube}).

\subsection{An upper bound for $\mu(Q_d)$}
Here, we first analyze specific optimal solutions for the \textsc{Mutual-Visibility} problem in hypercubes of small dimension $d\leq 5$ and then we provide an upper bound holding for any $d\ge 6$.

First of all, we can claim the following lemma.

\begin{lemma}\label{lemma:convex_hper}
Each subgraph $Q'$ of $Q_d$ that is isomorphic to $Q_{d'}$, $d'<d$, is convex. 
\end{lemma}

\begin{proof}
Consider a hypercube $Q'$ subgraph of $Q_d$ that is isomorphic to $Q_{d'}$ for some $d'<d$.
By considering two vertices $x$ and $y$ of $Q'$, we have to show that there not exists a shortest path between $x$ and $y$ passing through a vertex $z\in V(Q_d)\setminus V(Q')$.
Since the distance between two vertices in a hypercube is governed by the Hamming distance, the bits that differ between the labels associated with $x$ and $y$ concern only $Q'$. In fact, there cannot be a shortest path that makes use of a vertex $z\in V(Q_d)\setminus V(Q')$ as, by construction, the corresponding bit leading to (the dimension of) $z$ is different from the one in the labels of both $x$ and $y$.	\qed
\end{proof}

\begin{corollary}\label{cor:d}
    $\mu(Q_d) \le 2\mu(Q_{d-1})$, for each $d\ge 2$.	
\end{corollary}
\begin{proof}
    The proof simply follows by recalling that $Q_d$ can be obtained by the Cartesian product of $Q_{d-1} \cp K_2$, for $d\ge 2$, i.e., by suitably connecting two hypercubes of dimension $d-1$. Therefore, by Lemma~\ref{lemma:convex_hper} the claim holds. \qed
\end{proof}

We now consider all the dimensions $d\leq 5$, one by one, for which we can provide optimal solutions.

From results provided in~\cite{DiStefano22}, we know that $\mu(Q_1)=\mu(K_2) = 2$ and that  $\mu(Q_2)=\mu(C_4) = 3$.
For $Q_3$ (that contains $8$ vertices) there exists exactly one $\mu$-set of size 5, up to isomorphisms. By referring to Fig.~\ref{fig:hypercube}, the optimal solution is provided by the set $\{000$, $001$, $100$, $110$, $011\}$. 
It can be observed that there not exists any solution with more than 5 vertices. By contradiction, assume that $X$, with $|X|\ge 6$, is a $\mu$-set for $Q_3$. Then, for each subcube $Q'$ of $Q_3$ isomorphic to $Q_2$ we have $|V(Q')\cap X| < 4$ (in fact, $Q'$ is a convex subgraph isomorphic to a cycle). Hence, the two copies of $Q_2$ must have 3 elements each in $X$. 
To avoid having a cycle $C_4$ whose elements are all in $X$, the two vertices not in $X$ must be antipodal; but, this implies that there exists another pair of antipodal vertices that are in $X$ which are not mutually-visible.

\begin{figure}[h]
\centering
\includegraphics[scale=0.65]{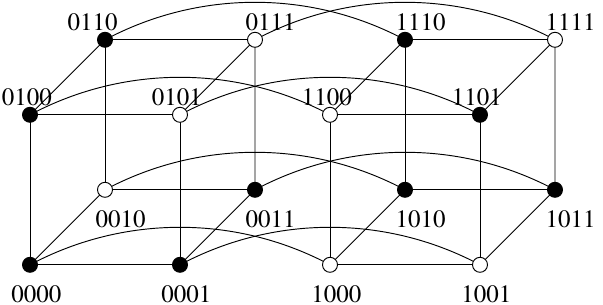}~~  
\caption{A representation of $Q_4$ along with a $\mu$-set for that graph (the black vertices).}
\label{fig:Q4}
\end{figure}

For $Q_4$ (that contains $16$ vertices), Corollary~\ref{cor:d} implies that each $\mu$-set can have at most 10 elements, i.e., 5 vertices per each subcube of dimension 3. Since the $\mu$-set for $Q_3$ is unique up to isomorphisms, we have few possibilities to combine the two subcubes. 
By a computer-assisted exhaustive search, we have that there not exists a solution with 10 vertices selected, whereas we can provide a mutual-visibility set with 9 vertices selected. By referring to Fig.~\ref{fig:Q4}, the optimal solution is provided by the set $\{0000$, $0001$, $0100$, $0110$, $0011$, $1101$, $1010$, $1011$, $1110\}$. 
Furthermore, the obtained solution is unique up to isomorphisms.

\begin{figure}[h]
\centering
\includegraphics[scale=0.65]{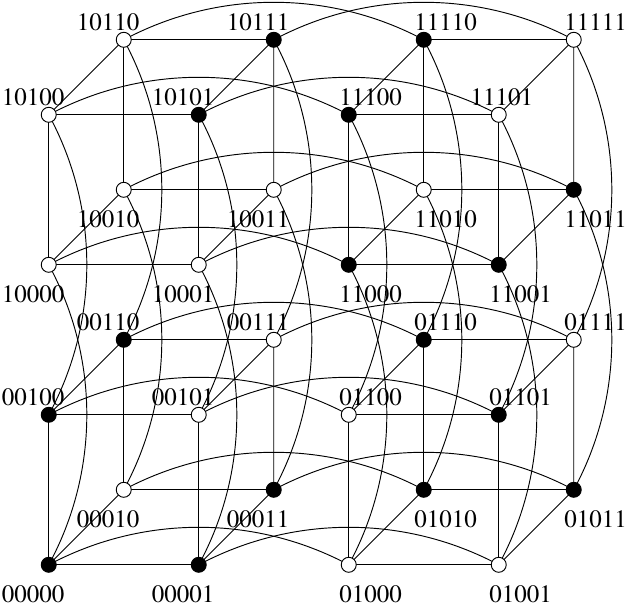}~~  
\caption{A representation of $Q_5$ along with a $\mu$-set for that graph (the black vertices).} 
\label{fig:Q5}
\end{figure}

Concerning $Q_5$, since $\mu(Q_4)=9$, by Lemma~\ref{lemma:convex_hper} we can obtain a mutual-visibility set with at most 18 vertices. 
Again, by means of a computer-assisted case-by-case analysis,  
we have that the optimal and unique (up to isomorphisms) solution is of size 16. By referring to Fig.~\ref{fig:Q5}, the optimal solution is provided by the set $\{00000$, $00001$, $00100$, $00110$, $00011$, $01101$, $01010$, $01011$, $01110$, $10101$, $10111$, $11000$, $11001$, $11100$, $11110$, $11011\}$.

For $d\ge 6$, we can state the next corollary that provides an upper bound to $\mu(Q_d)$.

\begin{corollary}\label{cor:upbound}
%
$\mu(Q_d) \le 2^{d-1}$, for each $d\ge 5$.	
\end{corollary}

\begin{proof}
We have shown that for $Q_5$ any $\mu$-set contains $2^4$ elements, i.e., exactly half of the vertices. The claim then simply follows by Corollary~\ref{cor:d}.
\qed
\end{proof}

Summarizing, in Table~\ref{tab:optimal} we report all the obtained results concerning the upper bounds for $\mu(Q_d)$.

\begin{table}
\begin{center}
\setlength{\tabcolsep}{6pt}
\renewcommand{\arraystretch}{1.2}
\begin{tabular}{ |c|c|l| } 
 \hline
 $d$ & $n(Q_d)$ & $\mu(Q_d)$  \\ 
 \hline
 1 & 2  & $\mu(Q_d)=2$ \\ 
 2 & 4  & $\mu(Q_d)=3$ \\
 3 & 8  & $\mu(Q_d)=5$ \\
 4 & 16 & $\mu(Q_d)=9$ \\ 
 5 & 32 & $\mu(Q_d)=16$ \\
 $d\ge 6$ & $2^{d}$ & $\mu(Q_d) \le 2^{d-1}$ \\
 \hline
\end{tabular}
\end{center}\caption{On the size of a mutual-visibility set $X$ for $Q_d$ as $d$ varies. We found $\mu$-sets (i.e., optimal solutions) for $d\leq 5$ and an upper bound for $\mu(Q_d)$ when $d\ge 6$. Note that the $\mu$-set for $Q_d$, $d<6$, is unique up to isomorphisms.}
\label{tab:optimal}
\end{table}

\subsection{Lower bound and an approximation algorithm}

Here, we first provide a lower bound for $Q_d$ and then we derive an approximation algorithm for the mutual-visibility problem in the class of hypercubes.

\begin{theorem}\label{thm:lb-Qd}
$\mu(Q_d)\ge 
\binom{d}{\lfloor\frac d 2 \rfloor} + 
\binom{d}{\lfloor\frac d 2 \rfloor+3}$, for any $d\ge 1$.
\end{theorem}
\begin{proof}
Let $X_p$ be the subset of $V(Q_d)$ containing all the vertices that are at distance $p$ from the vertex labeled with all zeroes. By construction, and by the Hamming distance, the elements of $X_p$ are all the vertices whose labels contain exactly $p$ 1's.
Hence, it is easy to find a shortest path between two of the selected vertices, say $x$ and $y$ that are at a distance $j$ from each other. It suffices to detect the $j$ differences among the labels associated with the two vertices. Then, by first replacing (one per step) the 1's present in $x$ but not in $y$ with 0's and then similarly replacing the 0's with 1's, equals to determine a shortest path between $x$ and $y$. Such a path is shortest because at each step makes a change in the direction of the destination, i.e., its length is exactly $j$. Furthermore, along the chosen path, there are no vertices in $X_p$ as the number of 1's is always smaller than $p$ until the destination.

The number of vertices with a fixed number p of 1's is $\binom{d}{p}$, which is maximized for $p=\lfloor\frac d 2 \rfloor$.

For $d < 5$,   the size of $X_{\lfloor\frac d 2 \rfloor}$ is $\binom{d}{\lfloor\frac d 2 \rfloor}$ and assumes the values $1,2,3,6$ for $d=1,2,3,4$. These values are less than the corresponding values of $\mu(Q_d)$ shown in Table~\ref{tab:optimal}. 
Since in these cases $\binom{d}{\lfloor\frac d 2 \rfloor+3}=0$, the statement holds.

For $d\ge 5$, we can choose a larger set of vertices, still guaranteeing mutual visibility. Given a hypercube $Q_d$, we define a set $X $ as follows:

\begin{itemize}
\item Let $p=\lfloor\frac d 2 \rfloor$;
\item $X = X_p\cup X_{p+3}$.
\end{itemize}

Let $x,y\in X$. There are shortest path between $x$ and $y$ whose internal vertices are labeled with $p+1$ and $p+2$ 1's only. Indeed, a shortest path from a vertex labeled with $p$ 1's to one labeled with $p+3$ 1's can involve only internal vertices obtained by first changing two 0's in 1's and then alternating one 1 in 0 and one 0 in 1 until the last step. By the number of 1's in each label of the internal vertices we are guaranteed that $x$ and $y$ are mutually visible. If  $x,y\in X_p$, to show the mutual visibility we can find a shortest path with vertices not in $x$ as describe above, whereas if $x,y\in X_{p+3}$ the vertices not in $X$ in a shortest path between $x$ and $y$, can be obtained by first replacing (one per step) the 0’s present in $x$ but not in $y$ with 1’s and then similarly replacing the 1’s with 0’s.
\qed
\end{proof}

An immediate consequence of the above theorem is the existance of an approximation algorithm for the mutual-visibility problem in the context of hypercubes, as stated by the following two corollaries.












\begin{corollary}\label{cor:approx-Qd}
There exists an algorithm for the \textsc{Mutual-Visibility} problem on a hypercube $Q_d$ which provides a solution of size greater than $\frac {2^{d}} {\sqrt{\frac \pi 2d}}$.
\end{corollary}
\begin{proof}
For $d\le 5$, we have already shown we can provide the optimal solution. For $d\ge 6$, 
we show that the procedure defined in the proof of Theorem~\ref{thm:lb-Qd} (for defining the $\mu$-set $X$ for $Q_d$) provides the requested algorithm. 

The cardinality of $X$ is equal to

$$|X|= \binom{d}{p} + \binom{d}{p+3}=
\binom{d}{\lfloor\frac d 2 \rfloor} + \binom{d}{\lfloor\frac d 2 \rfloor+3}.$$

By the Stirling's approximation~\cite{stirling1749differential}, we have that, for \emph{sufficiently large} values of $d$, 

$$\binom{d}{\lfloor\frac d 2 \rfloor} \sim \frac {2^{d}} {\sqrt{\frac \pi 2d}} \text{\ \ \ and consequently \ \ } |X| \sim \frac {2^{d+1}} {\sqrt{\frac \pi 2d}}.$$

However, for $d\ge 6$, we can easily verify that:

$$|X|> \frac {2^{d}} {\sqrt{\frac \pi 2d}},$$

since the inequality holds for $d=6$, and from there on the gap between the two terms monotonically increases, then the statement holds. 
\qed
\end{proof}

\begin{corollary}
    There exists an $O(\sqrt{d})$-approximation algorithm for the \textsc{Mutual-Visibility} problem on a hypercube $Q_d$.
\end{corollary}
\begin{proof}
By comparing the upper bound provided by Corollary~\ref{cor:upbound}  with the lower bound provided by Corollary~\ref{cor:approx-Qd}, we have an approximation ratio of:

$$\frac {\mu(Q_d)} {|X|} \leq \frac  {2^{d-1}\sqrt{\frac \pi 2 d}}{2^{d}}= O(\sqrt{d}).$$
\qed

\end{proof}

We remind that in a hypercube $Q_d$ of $n$ vertices, $n=2^d$, i.e., $d=\log n$. Consequently, the approximation provided by the above theorem for an $n$-vertex hypercube can be expressed as $O(\sqrt{\log n})$.


\subsection{Total mutual-visibility} 
In~\cite[Theorem 3]{BujtaKT23}, authors provide a lower bound for total mutual-visibility number of $K_{s}^{\square,r}$, where $K_{s}^{\square,r}$ denotes the Cartesian product of $r$ copies of $K_s$. Their proof is based on a probabilistic approach similar to the idea used in the proof in~\cite[Section 4]{BES73} for a famous hypergraph Turán-problem of Brown, Erdős, and Sós. The provided lower bound is  
\begin{equation}\label{eq:lb-sandi}
  \mu_t(K_{s}^{\square,r}) \ge \frac{s^{r-2}}{r(r+1)}. 
\end{equation}
Since for $s=2$ and $r=d$ the graph $K_{s}^{\square,r}$ corresponds to $Q_d$, the above result reduces to the following bound for hypercubes: 
$\mu_t(Q_d) \ge \frac{2^{d-2}}{d(d+1)}$.
It is worth remarking that this bound is asymptotically much smaller than the bound 
$\mu(Q_d) > \frac {2^{d}} {\sqrt{\frac \pi 2d}}$ provided in the proof of Corollary~\ref{cor:approx-Qd}. In~\cite{BujtaKT23}, the bound in Eq.~\ref{eq:lb-sandi} is used to prove that $\mu_t(K_{s}^{\square,r}) = \Theta(s^{r-2})$.

\section{Cube-connected cycles}\label{sec:ccc}
%

%

A cube-connected cycle of order $d\ge 3$ (denoted $\CCC_d$) can be defined as an undirected graph formed from a set of $d\cdot 2^d$ vertices labelled $[\ell,x]$, where $\ell$ is an integer between 0 and $d-1$ and $x$ is a binary string of length $d$. 

Two vertices $[\ell,x]$ and $[\ell',x']$ are adjacent if and only if: (i) either $x= x'$ and $|\ell-\ell'|=1$,%
\footnote{All arithmetic on indices and levels concerning $\CCC_d$ is assumed to be module $d$.}
 (ii) or $\ell = \ell'$ and $x' =x(\ell)$. In this last case, $x$ and $x'$ differ exactly for the bit in position $\ell$.



The edges which satisfy condition (i) are referred to as \emph{cycle edges}, whereas the edges that satisfy  condition (ii) are referred to as the \emph{hypercube edges}. 

Note that, the removal of all hypercube edges produces a graph with $2^d$ components, each of which is isomorphic to a cycle $C_d$. For this reason, each cycle $C_d$ in $\CCC_d$ which does not include any hypercube edge is referred to as a \emph{supervertex} of the (embedded) hypercube.

Furthermore, contracting all the cycle edges in $\CCC_d$ will produce a graph with $2^d$ vertices 
isomorphic to a hypercube $Q_d$. 
Consequently, the smallest $\CCC_d$ is in fact defined for $d=3$.




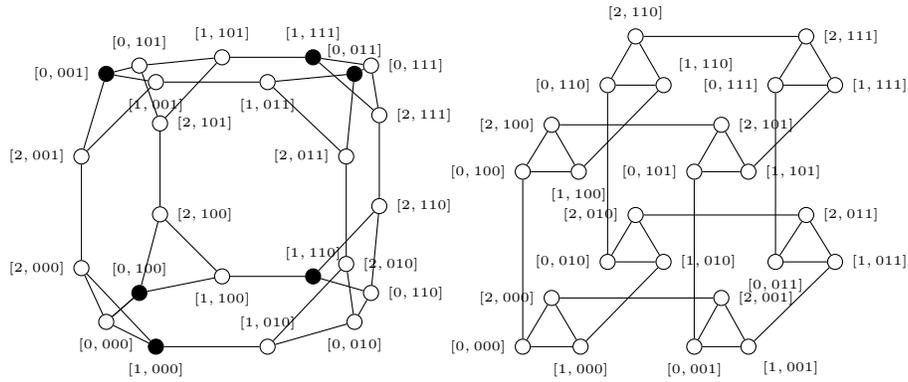
\begin{figure}
\begin{subfigure}{0.25\textwidth}

\begin{tikzpicture}[scale=0.11,font=\tiny]
\tikzstyle{nod}= [circle, draw,inner sep=0pt, minimum size=0.2cm]
\pgfsetxvec{\pgfpoint{1.cm}{0.0cm}}
\pgfsetyvec{\pgfpoint{0.0cm}{1.0cm}}
\node[nod] at (0,0) [label=below:${[0,000]}$] (a) {};
\node[nod] at (4,3.5) [fill=black, label=above:${[0,100]}$] (b) {};
\node[nod] at (14,5.5) [label=below:${[1,100]}$] (c) {};
\node[nod] at (25,5.5) [fill=black,label=above:${[1,110]}$] (d) {};
\node[nod] at (32,3.5) [label=right:${[0,110]}$] (e) {};
\node[nod] at (30,0) [label=below:${[0,010]}$] (f) {};
\node[nod] at (19.5,-3) [label=above:${[1,010]}$] (g) {};
\node[nod] at (6,-3) [fill=black,label=below:${[1,000]}$] (h) {};
\node[nod] at (-3,6.5) [label=left:${[2,000]}$] (i) {};
\node[nod] at (6.5,13) [label=right:${[2,100]}$] (l) {};
\node[nod] at (33,14) [label=right:${[2,110]}$] (m) {};
\node[nod] at (29,7) [label=right:${[2,010]}$] (n) {};
\node[nod] at (-3,20) [label=left:${[2,001]}$] (o) {};
\node[nod] at (6.5,24) [label=right:${[2,101]}$] (p) {};
\node[nod] at (33,25) [label=right:${[2,111]}$] (q) {};
\node[nod] at (29,20) [label=left:${[2,011]}$] (r) {};
\node[nod] at (0,30) [fill=black, label=left:${[0,001]}$] (t) {};
\node[nod] at (4,31) [label=above:${[0,101]}$] (u) {};
\node[nod] at (14,32) [label=above:${[1,101]}$] (v) {};
\node[nod] at (25,32) [fill=black,label=above:${[1,111]}$] (z) {};
\node[nod] at (32,31) [label=right:${[0,111]}$] (y) {};
\node[nod] at (30,30) [fill=black,label=above:${[0,011]}$] (w) {};
\node[nod] at (19.5,29) [label=below:${[1,011]}$] (x) {};
\node[nod] at (6,29) [label=below:${[1,001]}$] (j) {};

\path (a)
 edge (b)
edge (h)
edge(i)
(b)
edge (a)
edge (c)
edge (l)
(c)
edge (l)
edge (d)
(d)
edge (m)
edge (e)
(e)
edge (m)
edge (f)
(f)
edge (n)
edge (g)
(g)
edge (h)
edge (n)
(h)
edge (i)

(o)
edge(t)
edge (j)
edge (i)
(r)
edge(x)
edge(w)
edge(n)
(p)
edge(u)
edge(v)
edge(l)
(q)
edge(m)
edge(y)
edge (z)

(t)
edge(u)
edge(j)
(v) 
edge (z)
edge (u)
(z) edge (y)
(y) edge (w)
(w) edge (x)
(x) edge (j);
\end{tikzpicture}

    \label{fig:first}
\end{subfigure}
~~~~~~~~~~~~~~~~~~~~~~~~
\begin{subfigure}{0.3\textwidth}

\begin{tikzpicture}[scale=0.12,font=\tiny]
\tikzstyle{nod}= [circle, draw,inner sep=0pt, minimum size=0.2cm]
\pgfsetxvec{\pgfpoint{1.cm}{0.0cm}}
\pgfsetyvec{\pgfpoint{0.0cm}{1.0cm}}

\node[nod] at (0,0) [label=left:${[0,000]}$] (a) {};
\node[nod] at (6.4,0) [label=below:${[1,000]}$] (b) {};
\node[nod] at (3.2,5.4) [label=left:${[2,000]}$] (c) {};

\node[nod] at (3.2,24.6) [label=left:${[2,100]}$] (d) {};
\node[nod] at (6.2,19.4) [label=below:${[1,100]}$] (e) {};
\node[nod] at (0,19.4) [label=left:${[0,100]}$] (f) {};

\node[nod] at (12.5,14.6) [label=left:${[2,010]}$] (g) {};
\node[nod] at (15.6,9.4) [label=right:${[1,010]}$] (h) {};
\node[nod] at (9.4,9.4) [label=left:${[0,010]}$] (i) {};

\node[nod] at (31.4,14.6) [label=right:${[2,011]}$] (l) {};
\node[nod] at (34.6,9.4) [label=right:${[1,011]}$] (m) {};
\node[nod] at (28,9.4) [label=below:${[0,011]}$] (n) {};

\node[nod] at (22,5.4) [label=right:${[2,001]}$] (o) {};
\node[nod] at (25,0) [label=below right:${[1,001]}$] (p) {};
\node[nod] at (19,0) [label=below:${[0,001]}$] (q) {};

\node[nod] at (22,24.6) [label=right:${[2,101]}$] (r) {};
\node[nod] at (25,19.4) [label=right:${[1,101]}$] (s) {};
\node[nod] at (19,19.4) [label=left:${[0,101]}$] (t) {};

\node[nod] at (12.5,34.4) [label=above:${[2,110]}$] (u) {};
\node[nod] at (15.6,29) [label=above right:${[1,110]}$] (v) {};
\node[nod] at (9.4,29) [ label=left:${[0,110]}$] (z) {};

\node[nod] at (31.4,34.4) [label=right:${[2,111]}$] (x) {};
\node[nod] at (34.6,29) [label=right:${[1,111]}$] (k) {};
\node[nod] at (28,29) [label=left:${[0,111]}$] (j) {};

\path (a)
edge (b)
edge(c)

(c)
edge (b);

\path (d)
edge (e)
edge(f)

(f)
edge (e)

(a)
edge (f);

\path (g)
edge (h)
edge(i)

(h)
edge (i)
edge (b);

\path (l)
edge (m)
edge(n)
edge(g)

(m)
edge (n);

\path (o)
edge (p)
edge(q)
(p)
edge (q)

(c)
edge (o)

(p)
edge (m);

\path (r)
edge (s)
edge(t)
(s)
edge(t)

(r)
edge (d)

(t)
edge (q);

\path (u)
edge (v)
edge(z)
(v)
edge (z)
(v)
edge (e)
(z)
edge (i);

\path (x)
edge (k)
edge(j)
(k)
edge(j)
(x)
edge(u)
(k)
edge(s)
(j)
edge(n);
\end{tikzpicture}

\label{fig:second}
\end{subfigure}        
\caption{\emph{(left)} A representation of a $CCC_3$ graph, black vertices represent a $\mu$-set for $CCC_3$. \emph{(right)} An alternative representation of a $CCC_3$. 
}
\label{fig:ccc}
\end{figure}

\paragraph{Natural routing in $CCC_d$.} 
Let $u=[\ell,x]$ and $v=[\ell',x']$ be two vertices of a $CCC_d$, and $H=\{p_0,\ldots, p_{h-1}\}$ be the positions in which the binary strings representing $u$ and $v$ differ. By~\cite{V93}, the distance between $u$ and $v$ equals $h+k$, with $k$ being the number of edges in the shortest walk on a cycle $C_d$ starting at index $\ell$ and ending at index $\ell'$ which includes a visit to every vertex with an index in the set $H$. Then, a \emph{natural routing} from $u$ to $v$ is to traverse the shortest walk while traversing the hypercube edges each time a vertex with an index in the set $H$ is met for the first time. 

\begin{lemma}\label{lemma:convex_ccc}
Each subgraph $Q'$ composed by the supervertices of a $CCC_d$, isomorphic to $Q_{d'}$, $d'<d$, induces a subgraph of $CCC_d$ that is convex. 
\end{lemma}

\begin{proof}
By Lemma~\ref{lem:convex-subgraph}, we have that the shortest paths between vertices in $Q_{d'}$ are confined within $Q_{d'}$ itself.
Consider now the subgraph obtained by restoring each supervertex with  the corresponding cycle $C_d$ along with its edges within $Q_{d'}$. By the natural routing, any shortest path among two vertices $u$ and $v$ of the defined subgraph makes use of the hypercube edges only from vertices with index in the set $H$. Since all such edges are included in the obtained subgraph, the claim holds.
	\qed
\end{proof}

\begin{lemma}\label{lemma:2upbound-ccc}
$\mu(\CCC_d) \le 3\cdot 2^{d-2}$.	
\end{lemma}

\begin{proof}
The proof follows by first observing that for $d=3$, the optimal solution provides $\mu(CCC_3)=3\cdot 2^{3-2}=6$. In fact, we have obtained this result by a computer-assisted exhaustive search, and the optimal solution found is shown in Fig.~\ref{fig:ccc}. 

Then, by Lemma~\ref{lemma:convex_ccc}, we have that, as $d$ increases, $\mu(CCC_d)$ can double at most.\qed
\end{proof}

\begin{lemma}\label{lemma:lb-ccc}
$\mu(CCC_d)\ge 
2^{\lceil \frac d 2 \rceil-1}$, for any $d\ge 3$.
\end{lemma}
\begin{proof}
Given a cube connected cycle $CCC_d$, $d \ge 3$, we define a set $X$ according to the following procedure:

\begin{itemize}
\item Insert into $X$ any vertex $v=[0,x]$ with a number of 1's bounded by $\lceil \frac d 2 \rceil$ in the least significative positions, but position 0, of $x$.
\end{itemize}


In doing so, and reminding the natural routing in a $CCC_d$, we have that the shortest path between two vertices $v_1$ and $v_2$ belonging to $X$ is well defined by considering the walk in $C_d$ that from $0$ reaches $p_{h-1}$ and then comes back to $0$ either proceeding backward or forward, it depends on the distances. Since all 
vertices in $X$ have their first coordinate 0, all the traversed edges within any encountered $C_d$ as well as all hypercube edges included in the shortest path never meet other vertices in $X$. 

Summarizing, we have that any two vertices in $X$ are mutually visible since the shortest path forced by the chosen set $X$ ensures to never encounter a vertex $[0,z]$ for any $z\not\in \{x,y\}$.
By construction, $|X|= 2^{\lceil \frac d 2 \rceil -1}$. \qed
\end{proof}

\begin{theorem}\label{th:approx-ccc}
There exists a $3 \cdot 2^{\lfloor \frac d 2 \rfloor-1}$-approximation algorithm for the \textsc{Mutual-Visibility} problem on a cube-connected cycle $CCC_d$.
\end{theorem}

\begin{proof}

By Lemmata~\ref{lemma:2upbound-ccc} and~\ref{lemma:lb-ccc}, we obtain the claimed approximation ratio for the algorithm provided in the proof of Lemma~\ref{lemma:lb-ccc}:
$$ \frac{\mu(\CCC_d)}{|X'|} \leq 
   \frac {3\cdot 2^{d-2}} { 2^{\lceil \frac d 2 \rceil-1}}= 3 \cdot 2^{\lfloor \frac d 2 \rfloor-1}.
$$
\qed
\end{proof}

We remind that in a cube-connected cycle $CCC_d$ of $n$ vertices, $n=d\cdot 2^d$, i.e., $d=\log n - \log d = \log n -o(\log n)$. Consequently, the approximation provided by the above theorem for an $n$-vertex cube-connected cycle can be expressed as $O(\sqrt{n})$.

\subsection{Total mutual visibility}
In~\cite{totalmutualzero}, a characterization of all graphs having total mutual-visibility zero is provided. The characterization is based on the notion of \emph{bypass vertex}: a vertex $u$ of $G$ is a bypass vertex if $u$ is not the middle vertex of a convex $P_3$ of $G$. Let $\bp(G)$ be the number of bypass vertices of $G$. 

\begin{theorem}\label{cor:total-ccc}
\emph{~\cite[Theorem 3.3 ]{totalmutualzero}}
Let $G$ be a graph with $n(G)\ge  2$. Then $\mu_t(G) = 0$ if and only if $\bp(G) = 0$.
\end{theorem}

This result directly provides the value of $\mu_t(\CCC_d)$.

\begin{corollary}
$\mu_t(\CCC_d)=0$, for each $d\ge 3$.
\end{corollary}
\begin{proof}
Let $u$ be a vertex of $\CCC_d$ and consider a path $P_3$ given by $(v',u,v'')$, where $v'$ and $v''$ belong to different cycles of the graph. It can be observed that such a path is convex. Since it is well-known that $\CCC_d$ is a vertex-transitive graph, then each vertex $u$ is not a bypass vertex and hence $\bp(CCC_d)=0$.
\end{proof}

\section{Butterfly}\label{sec:butterflies}
\begin{figure}[h]
\begin{center}
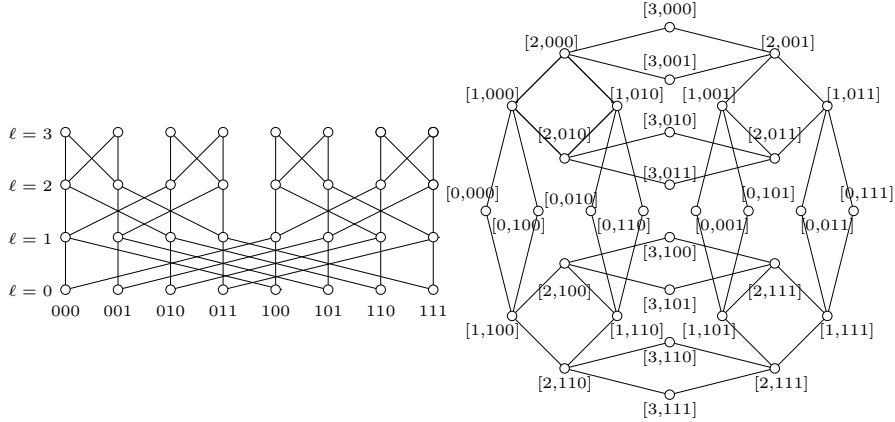
\caption{Normal representation and diamond representation of $\BF(3)$.}
\label{fig:butterfly}
\end{center}
\end{figure}

%
A $d$-dimensional butterfly $\BF(d)$ is an undirected graph with vertices $[\ell,c]$, where $\ell\in \{0,1,\ldots,d\}$ is the \emph{level} and $c\in \{0,1\}^d$ is the \emph{column}. The vertices $[\ell,c]$ and $[\ell',c']$ are adjacent if $|\ell - \ell'| = 1$, and either $c = c'$ or $c$ and $c'$ differ precisely in the $\ell$-th bit. $\BF(d)$ has $d + 1$ levels with $2^d$ vertices at each level, which gives $(d+1)\cdot 2^d$ vertices in total. The vertices at level $0$ and $d$ are $2$-degree vertices and the rest are $4$-degree vertices. $\BF(d)$ has two standard graphical representations, namely normal and diamond representations (see Fig.~\ref{fig:butterfly}). For further details one may refer to~\cite{ManuelARR08}. 

\begin{figure}[h]
\centering
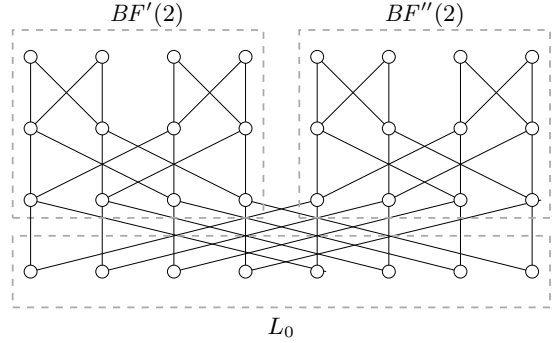
\caption{Partition of $\BF(3)$ into $\BF'(2)$, $\BF''(2)$, and $L_0$.}
\label{fig:bf-partition}
\end{figure}

Some additional notation is required: let $A_i = \{ [\ell,i]~|~ \ell\in \{0,1,\ldots,d\} \}$ be the vertex set forming the $i$-th column of the butterfly, and let $L_i = \{ [j,c]~|~ c\in \{0,1\}^d \}$, be the vertex set forming its $j$-th level. Note that, $\BF(d)$ can be partitioned into two copies of $\BF(d-1)$ (that we denote as $\BF'(d-1)$ and $\BF''(d-1)$) and $L_0$ (cf. Fig.~\ref{fig:bf-partition}). Notice that both $\BF'(d-1)$ and $\BF''(d-1)$ are convex subgraphs of $\BF(d)$.

\paragraph{Natural routing in $\BF(d)$.} Consider the case in which starting from a vertex $[0,i]$, it is necessary to reach a vertex $[d,j]$ along a shortest path. A \emph{natural routing} for such a task simply requires comparing the corresponding bits of $i$ and $j$ starting from the leftmost: if they coincide, the edge along the current column is traversed, otherwise, the edge for changing column is used. Symmetrically, a routing from $[d,i]$ toward $[0,j]$ can be obtained in a similar way but comparing the corresponding bits of $i$ and $j$ starting from the rightmost. This leads to the following useful property:
\begin{description}
\item[$\Prop$:]
Let $u=[\ell',i]$ and $v=[\ell'',j]$ be two vertices, $\ell_{\min} = \min\{\ell',\ell''\}$ and $\ell_{\max} = \max\{\ell',\ell''\}$. 
Each shortest path between $u$ and $v$ is comprised between levels $\min \{k'-1,\ell_{\min} \}$ and $\max \{k'',\ell_{\max} \}$, where $k'$ and $k''$ are the positions of the first and last bit, respectively, in which $i$ and $j$ differ.
\end{description}
As special cases for this property we get:
($\Prop_1$) if from $[0,i]$ it is necessary to reach $[0,j]$, then a shortest path reaches level $k$ iff the $k$-th bits of $i$ and $j$ differ; ($\Prop_2$) if from $[d,i]$ it is necessary to reach $[d,j]$, then a shortest path reaches level $d-k$ iff the $(d-k)$-th bits of $i$ and $j$ differ.

\begin{lemma}\label{lem:at-most-two-per-column}
Let $X$ be a mutual-visibility set of $\BF(d)$. Then,  $|X\cap A_i| \le 2$ for each $i\in \{0,1,\ldots,d\}$.
\end{lemma}
\begin{proof}
According to Lemma~\ref{lem:convex-subgraph}, since $A_i$ is a convex subgraph, $A_i$ is isomorphic to a path graph, and the mutual-visibility number of a path graph is two.
\qed
\end{proof}

\begin{lemma}\label{lem:mutual-visibility-set-for-BF}
$X= (L_0\cup L_d) \setminus \{ [0,111\ldots 1], [d,111\ldots 1]\} $
is a mutual-visibility set of $\BF(d)$.
\end{lemma}
\begin{proof}
Consider $u=[0,i]$ and $v=[d,j]$. A shortest path that makes $u$ and $v$ $X$-visible is given by the natural routing from $u$ to $v$. 

Consider now $u=[0,i]$ and $v=[0,j]$. According to Property $\Prop_1$, if the last bit of $i$ and $j$ does not differ, then there exists a shortest path whose interior vertices belong all at levels $1,2,\ldots, d-1$, and this property trivially makes the two vertices $X$-visible. If the last bit of $i$ and $j$ differs, then their distance is $2d$ and a shortest path that makes the two vertices $X$-visible can be defined by composing the natural routing from $u$ to $[d,111\ldots 1]$ with the natural routing from $v$ to $[d,111\ldots 1]$. 

Consider the last case in which $u=[d,i]$ and $v=[d,j]$. According to Property $\Prop_2$, if the first bit of $i$ and $j$ does not differ, then there exists a shortest path whose interior vertices belong all at levels $1,2,\ldots, d-1$, and this property trivially makes the two vertices $X$-visible. If the first bit of $i$ and $j$ differs, then their distance is $2d$ and a shortest path that makes the two vertices $X$-visible can be defined by composing the natural routing from $[0,111\ldots 1]$ to $u$ with the natural routing from $[0,111\ldots 1]$ to $v$.
\qed
\end{proof}

\begin{lemma}\label{lem:upper-bound-for-BF}
$\mu( \BF(d) ) \le 2^{d+1}-2$.
\end{lemma}
\begin{proof}
By contradiction, assume the statement is false, and let $X$ be a mutual-visibility set of $\BF(d)$. We can restrict the analysis to the case $|X|= 2^{d+1}-1$ because, if the size is $2^{d+1}$, then the removal of one element from $X$ would still produce a mutual-visibility set. 
%
%
According to the cardinality of $X$ and Lemma~\ref{lem:at-most-two-per-column}, we get %
\begin{equation}\label{eq:one-at-tau}
    |X\cap A_{\tau}| = 1 \textrm{ for some } \tau\in \{0,1\}^d,
\end{equation}
and
\begin{equation}\label{eq:two-at-remaining}
    |X\cap A_i| = 2 \textrm{ for each } i\in \{0,1\}^d \setminus \{\tau\}.
\end{equation}
If $A_i$, for each $i\in \{0,1\}^d$, shares two elements with $X$, we denote them as $[X_{\max},i]$ and $[X_{\min},i]$, where $[X_{\max},i]$ ($[X_{\min},i]$, resp.) is the element in $|X\cap A_i|$ at the largest (smaller, resp.) level. If $A_i$ shares just one element with $X$, then $[X_{\min},i] = [X_{\max},i]$.

\medskip
\textbf{Case 1: $|X\cap L_0| \ge 2$.} 
Let $A_i$ be a column such that $[X_{\min},i] = [0,i]$, and let $[1,j]$ be the neighbour of $[0,i]$ belonging $\BF''(d-1)$. Notice that $|X\cap L_0| \ge 2$ allow us to assume  $[X_{\min},j] \neq [X_{\max},j]$ (in fact, if $j=\tau$ we can start this analysis from another column $A_{i'}$ such that $[X_{\min},i'] = [0,i']$). Without loss of generality, assume that column $A_i$ intersects $\BF'(d-1)$. 

It can be easily observed that if the level of $[X_{\min},j]$ is greater or equal to 1, then  $[X_{\min},i]$ and $[X_{\max},j]$ are not $X$-visible since they are obstructed by $[X_{\min},j]$. As a consequence, we get that $[X_{\min},j] = [0,j]$. By repeating the previous arguments, from $[X_{\min},j] = [0,j]$ we get that there exists $[X_{\min},k]$, for some $k\not\in \{i,j\}$, that is located at level 0. Iteratively, we get that 

\begin{equation}\label{eq:all-but-one-at-zero}
  [X_{\min},i] = [0,i] \textrm{ for each } i\in \{0,1\}^d\setminus \{\tau\}.
\end{equation}

If also $[X_{\min},\tau] = [0,\tau]$ holds, then the entire level 0 of $\BF(d)$ belongs to $X$. In this case, by selecting a vertex $[X_{\max},i']$ from $\BF'(d-1)$ and a vertex $[X_{\max},i'']$ from $\BF''(d-1)$ we easily get that each shortest path between them passes through a vertex at level 0. This implies that $[X_{\max},i']$ and  $[X_{\max},i'']$ are not $X$-visible.

If $[X_{\min},\tau] \neq [0,\tau]$, then, without loss of generality, we assume that $A_{\tau}$ intersects $\BF''(d-1)$. In this case, in order to guarantee the mutual-visibility property of $X$, each shortest path from a vertex of $X$ belonging to $\BF''(d-1)$ to a vertex of $X$ belonging to $\BF'(d-1)$ must pass through the edge $([0,\tau], [1,t])$, where $A_t$ is a column intersecting $\BF'(d-1)$. In order to reach from such an edge each element of $X$ located in $\BF'(d-1)$, all such elements to be reached must necessarily be at level $d$. 

Observe that $[X_{\min},\tau] \neq [0,\tau]$ implies the membership of $[X_{\min},\tau]$ to $\BF''(d-1)$. 
Since $\BF''(d-1)$ is a convex subgraph, it can be noticed that in order to guarantee the mutual-visibility property of $X$ all the vertices of $X$ located in $\BF''(d-1)$ must be at level $d$. Hence we get  

\begin{equation}\label{eq:all-at-d}
  [X_{\max},i] = [d,i] \textrm{ for each } i\in \{0,1\}^d.
\end{equation}

According to Eq.~\ref{eq:all-at-d}, 
select now two elements of $X$ that are both at level 0 and at adjacent columns. These vertices have distance $2d$ and each shortest path connecting them passes through a vertex at level $d$. Since by Eq.~\ref{eq:all-at-d} each vertex at level $d$ is in $X$, then we get that the two selected vertices at level 0 are not $X$-visible.

\smallskip
\textbf{Case 2: $|X\cap L_0| \le 1$.} 
Let $X'$ ($X''$, respectively) be the set containing all vertices of $X$ that are contained in $\BF'(d-1)$ ($\BF''(d-1)$, respectively).  According to Lemma~\ref{lem:convex-subgraph} and since both $\BF'(d-1)$ and $\BF''(d-1)$ are convex subgraphs, then $X'$ is a mutual-visibility set of $\BF'(d-1)$ and $X''$ is a mutual-visibility set of $\BF''(d-1)$. Since the cardinality of both $X'$ and $X''$ is at least $2^d-1$, then we can recursively apply this proof to both the subgraphs. Then, either Case 1 will eventually apply, or we end the recursion with the terminal situation given by a subgraph isomorphic to $\BF(2)$ that shares with the original set $X$ seven or eight elements, at most two of such elements per column, and at most one element at the lower level. 

If $\BF(2)$ shares 8 vertices with $X$, it is easy to observe that the presence of at most one vertex of $X$ at the lower level implies the existence of a subgraph isomorphic to $\BF(1)$ with its 4 vertices all belonging to $X$. Since $\BF(1)$ is isomorphic to a cycle $C_4$ and $\mu(C_4)=3$, then we get a contradiction for Lemma~\ref{lem:convex-subgraph}.

\begin{figure}[h]
\centering
\includegraphics[width=9.5cm]{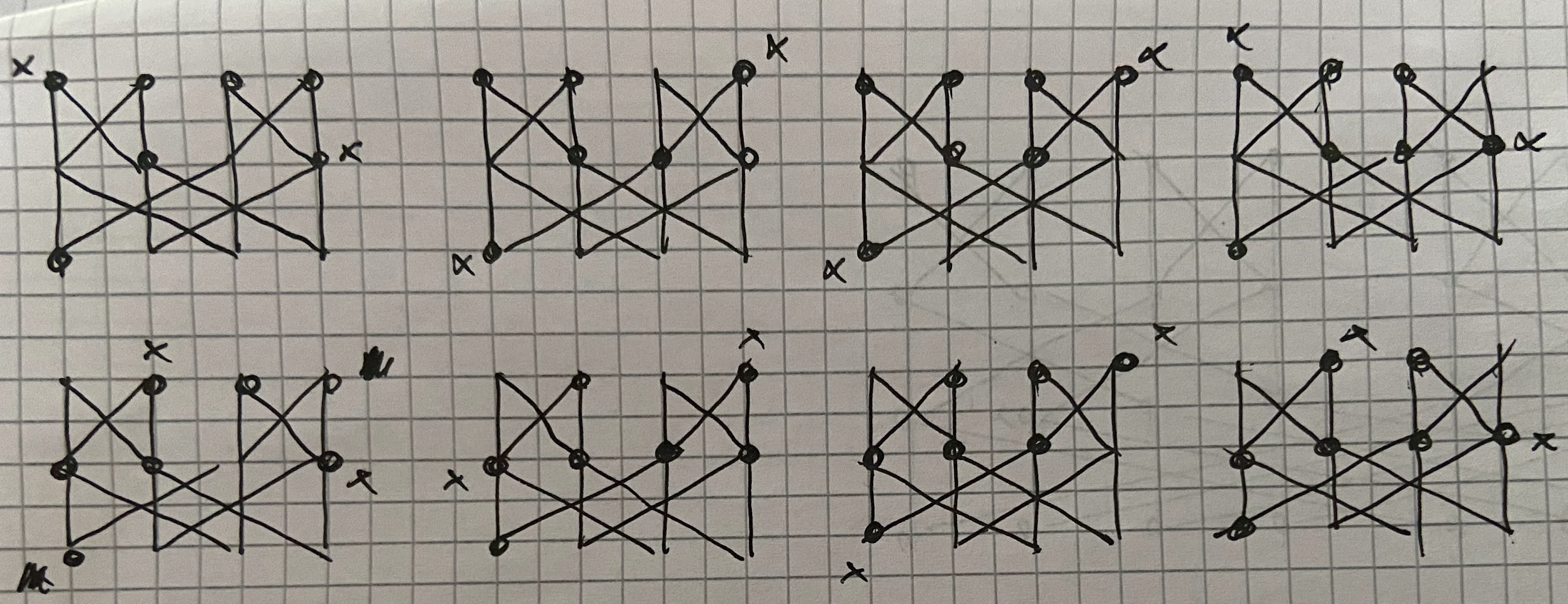}
\caption{Visualization of all the butterflies $\BF(2)$, up to isomorphisms, defined in the proof of Lemma~\ref{lem:upper-bound-for-BF}. For each graph, black vertices form a subset to be tested against mutual visibility, and the pair of vertices marked with `\textsf{x}' are not mutually visible.}
\label{fig:BF2}
\end{figure}

If $\BF(2)$ shares 7 vertices with $X$, then there are only few possible configurations for $\BF(2)$ that do not contain a subgraph $\BF(1)$ with its 4 vertices all in $X$. Fig.~\ref{fig:BF2} shows all the configurations that, up to isomorphisms, match all such conditions. In each case, we get that there are 2 vertices that are not in mutual visibility. Again,  this contradicts Lemma~\ref{lem:convex-subgraph}.
\qed
\end{proof}

\begin{theorem}
\label{thm:optimal-for-BF}
$\mu( \BF(d) ) = 2^{d+1}-2$, for each $d\ge2$.
\end{theorem}
\begin{proof}
Since the mutual-visibility set provided by Lemma~\ref{lem:mutual-visibility-set-for-BF} contains $2^{d+1}-2$ elements, then the statement directly follows from Lemma~\ref{lem:upper-bound-for-BF}.
\qed
\end{proof}

\subsection{Total mutual-visibility}

\begin{figure}[h]
\centering
\includegraphics[width=6cm]{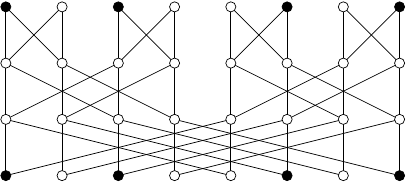}
\caption{Visualization of $L'_0\cup L'_d$ as computed in $\BF(3)$.
}
\label{fig:L0Ld}
\end{figure}

Given $\BF(d)$, let $L'_0$ ($L'_d$, respectively) be the subset of $L_0$ ($L_d$, respectively) containing all vertices $[0,i]$ ($[d,i]$, respectively) fulfilling one of the following conditions, where $i$ is interpreted as a binary number: 
\begin{itemize}
    \item $i$ is odd and $i\le 2^{d-1}$,
    \item $i$ is even and $i> 2^{d-1}$.
\end{itemize}
Fig.~\ref{fig:L0Ld} shows $L'_0\cup L'_d$ in $\BF(3)$.

\begin{lemma}\label{lem:total-mutual-visibility-set-for-BF}
$X = L'_0 \cup L'_d$ is a total mutual-visibility set of $\BF(d)$, for each $d\ge 1$.
\end{lemma}
\begin{proof}
let $u$ and $v$ be two distinct vertices of $\BF(d)$. From the proof of Lemma~\ref{lem:mutual-visibility-set-for-BF} we easily derive that $u$ and $v$ are $X$-visible in each of the following three cases: (1) $u=[0,i]$ and $v=[0,j]$, (2) $u=[0,i]$ and $v=[d,j]$, and (3) $u=[d,i]$ and $v=[d,j]$.

In order to prove the statement, it remains to consider three additional situations. 
\begin{itemize}
    \item[$(4)$]
    $u=[0,i]$ and $v=[\ell,j]$, with $\ell\not\in \{0,d\}$. According to the natural routing in $\BF(d)$, there exists a shortest path $P$ from $u=[0,i]$ to $u'=[d,j]$. Vertex $u=[\ell,j]$ belongs to $P$, and the $u,v$-subpath of $P$ guarantees that $u$ and $v$ are $X$-visible. 
    \item[$(5)$]
    $u=[d,i]$ and $v=[\ell,j]$, with $\ell\not\in \{0,d\}$. It is symmetric to the previous case. 
    \item[$(6)$]
    $u=[\ell',i]$ and $v=[\ell'',j]$, with $\ell',\ell''\not\in \{0,d\}$. 
    There are some subcases.
    \begin{itemize}
        \item[$(a)$]
        $i=j$. The statement trivially holds since a $u,v$-shortest path not passing through vertices of $X$ is contained in column $A_i$. 
        \item[$(b)$]
        $i\neq j$ and $u$ and $v$ are located on the same sub-butterfly. Without loss of generality, assume that both $u$ and $v$ are in $\BF'(d-1)$ and $i<j$.
        
        Since $u$ and $v$ are in the same sub-butterfly, each $u,v$-shortest path does not touch level 0 of $\BF(d)$. 
        
        If $i$ and $j$ are both even, then, according to Property $\Prop$, we get that each $u,v$-shortest path does not touch level $d$ of $\BF(d)$ unless $u$ or $v$ are at that level. Since $u$ and $v$ are not level $d$ by hypothesis, it follows that $u$ and $v$ are $X$-visible. When both $i$ and $j$ are odd their binary strings do not differ in the last position, and hence from Property $\Prop$ we still get that each $u,v$-shortest path does not touch level $d$; again, this guarantees the mutaul-visibility between the two vertices.

        If $i$ is odd and $j$ is even, then, according to Property $\Prop$, we get that each $u,v$-shortest path must touch level $d$ of $\BF(d)$, regardless of the level of $u$ an $v$; in this case, it can be observed that such a shortest path can pass through the vertex $[d,i]$. Conversely, if $i$ is even and $j$ is odd, the requested shortest path can pass through either $[d,i+1]$ or $[d,j]$. 
        \item[$(c)$]
        $i\neq j$, and $u$ and $v$ are located on different sub-butterflies. Without loss of generality, assume $u$ in $\BF'(d-1)$ and $v$ in $\BF''(d-1)$. If $i$ is odd, then it can be easily observed that a shortest path from $u$ to $v$ can be traversed by first reaching the vertex $[1,j]$ and then $[0,j']$ in the other sub-butterfly (note that this vertex is not in $X$ by definition). Now, the remaining part of the $u,v$-shortest path can be identified as in the case $(4)$ above. If $i$ is even, the requested shortest path can be defined by the same strategy, but now it is necessary to reach the vertex $[0,j]$ first and then $[1,j']$ in the other sub-butterfly. The remaining part of the $u,v$-shortest path can be identified as in the case $(6.a)$ or $(6.b)$.
    \end{itemize}
\end{itemize}
\qed
\end{proof}

\begin{theorem}\label{thm:total-optimal-for-BF}
$\mu_t( \BF(d) ) = 2^d$, for each $d\ge 1$.
\end{theorem}
\begin{proof}
Let $X$ be any total mutual-visibility set of $\BF(d)$. If $u=[0,i]$ and $v=[d,j]$, according to the natural routing of the butterfly we get that there exists a unique $u,v$-shortest path in $\BF(d)$. Consequently, each vertex in this $u,v$-shortest path cannot belong to $X$. In general, this leads to the property that each vertex $[\ell,i]$ such that $\ell\not\in \{0,d\}$ does not belong to $X$. 

Consider now two vertices $u=[d,i]$ and $v=[d,i+1]$, with $i$ even. They belong to a convex subgraph of the butterfly isomorphic to a cycle $C_4$, and hence they cannot both belong to $X$ otherwise the other pair of vertices in that cycle are not in mutual visibility. Symmetrically, vertices $u=[0,i]$ and $v=[0,i+2^{d-1}]$, $i$ even, belong to a convex subgraph $C_4$, and hence they cannot both belong to $X$.

From the arguments above it follows that in $X$ cannot exist any vertex with degree four. Moreover, only half of the vertices at level $d$ and at level 0 can belong to $X$. Hence, $\mu_t(\BF(d)) \le 2^d$. The statement follows by observing that the total mutual-visibility set provided by Lemma~\ref{lem:total-mutual-visibility-set-for-BF} contains exactly $2^d$ elements.
\qed
\end{proof}

\section{Conclusions}\label{sec:conclusions}
%
In this paper, we have studied the mutual-visibility and the total mutual-visibility in hypercubes, cube-connected cycles, and butterflies. While for any butterfly $\BF(d)$ we were able to provide exact formulae to calculate both $\mu(\BF(d))$ and $\mu_t(\BF(d))$, for the other topologies we were able to identify only approximation algorithms. These results, together with those obtained in~\cite{BujtaKT23} on Hamming graphs, suggest that the study of mutual visibility properties on such topologies seems to be particularly complex. This of course suggests further investigations within such topologies.


\end{document}